\documentclass[11pt,twoside,letter]{article}
\usepackage{amssymb,amsmath}
\usepackage{tikz}

\title{\bf\vspace{-39pt}Data driven sampling of oscillating signals}
\author{Brigitte Bidegaray-Fesquet \\
\small Univ. Grenoble Alpes, LJK, F-38000 Grenoble, France\\
\small CNRS, LJK, F-38000 Grenoble, France\\
\small Brigitte.Bidegaray@imag.fr \\
\\
Marianne Clausel \\
\small Univ. Grenoble Alpes, LJK, F-38000 Grenoble, France\\
\small CNRS, LJK, F-38000 Grenoble, France\\
\small Marianne.Clausel@imag.fr}
\date{}

\def\bbE{\mathbb{E}}
\def\bbN{\mathbb{N}}
\def\bbR{\mathbb{R}}
\def\bbZ{\mathbb{Z}}
\def\ph{\varphi}
\def\Osc{\operatorname{Osc}}
\newcommand\floor[1]{\left\lfloor #1 \right\rfloor}

\newtheorem{definition}{Definition}

\begin{document}

\maketitle

\begin{abstract}
The reduction of the number of samples is a key issue in signal processing for mobile applications. We investigate the link between the smoothness properties of a signal and the number of samples that can be obtained through a level crossing sampling procedure. The algorithm is analyzed and an upper bound of the number of samples is obtained in the worst case. The theoretical results are illustrated with applications to fractional Brownian motions and the Weierstrass function.
\vspace{5mm} 
\vspace{3mm}
\end{abstract}

\section{Introduction}

Autonomy, size and weight are very important issues in the design of mobile systems. One possibility to reduce the power consumption of signal processing systems is the reduction of the number of samples. Non uniform sampling is a way to have few samples for a large class of signals, especially sporadic signals, while still describing correctly the active parts of the signal. This leads to a smaller number of samples compared to Nyquist sampling \cite{Guan-Singer07, Mark-Todd81, Marvasti01, Sayiner-Sorensen-Viswanathan96}. Specific system architectures, such as event-driven architectures, allow the implementation of this specific sampling. These architectures take samples each time some specific event occurs, e.g.  specific voltage levels are crossed. Simple, low power, analog circuits can be designed to acquire information, possibly at high speed.

In this paper we want to relate the signal regularity to the number of non uniform samples obtained \textit{via} a level crossing technique. Indeed, intuitively, the more the signal is oscillating the more often the signal is sampled. This is of course a local property: the number of samples at the neighborhood of some point may then be related to the local smoothness of the signal, or more precisely to its H\"older regularity. This relationship will be tested on signals whose smoothness properties are perfectly known at each point. This can be useful to predict the processing complexity of e.g. biological signals such as EEG signals or fMRI data which are well-known to be both highly irregular and non stationary.

\begin{figure}[h]
\begin{center}
\begin{tabular}{cc}
\begin{tikzpicture}[scale=.6]
\draw[->,line width=1pt] (0,0) -- (10,0);
\draw[->,line width=1.5pt] (0,0) -- (0,6);
\draw[line width=1pt] (0,2) .. controls (2,5) and (3,7) .. (5,3)
            .. controls (5,3) and (8,-3) .. (10,3);
\foreach \x in {.5,1,1.5,2,2.5,3,3.5,4,4.5,5,5.5,6,6.5,7,7.5,8,8.5,9,9.5}
  \draw[dashed] (\x,0) -- +(0,6);
\draw[line width=1pt] (0,0) -- (10,0);
\draw[line width=1pt] (0,1) -- (10,1);
\draw[line width=1pt] (0,2) -- (10,2);
\draw[line width=1pt] (0,3) -- (10,3);
\draw[line width=1pt] (0,4) -- (10,4);
\draw[line width=1pt] (0,5) -- (10,5);
\draw(.5,2.75) circle (3pt);
\draw(1,3.5) circle (3pt);   \filldraw(1,3) circle (3pt);
\draw(1.5,4.2) circle (3pt); \filldraw(1.5,4) circle (3pt);
\draw(2,4.75) circle (3pt);
\draw(2.5,5.1) circle (3pt); \filldraw(2.5,5) circle (3pt);
\draw(3,5.2) circle (3pt);
\draw(3.5,5.05) circle (3pt);
\draw(4,4.6) circle (3pt);   \filldraw(4,5) circle (3pt);
\draw(4.5,3.9) circle (3pt); \filldraw(4.5,4) circle (3pt);
\draw(5,3) circle (3pt);     \filldraw(5,3) circle (3pt);
\draw(5.5,2.1) circle (3pt);
\draw(6,1.45) circle (3pt);  \filldraw(6,2) circle (3pt);
\draw(6.5,.95) circle (3pt); \filldraw(6.5,1) circle (3pt);
\draw(7,.6) circle (3pt);
\draw(7.5,.35) circle (3pt);
\draw(8,.35) circle (3pt);
\draw(8.5,.5) circle (3pt);
\draw(9,.95) circle (3pt);
\draw(9.5,1.8) circle (3pt); \filldraw(9.5,1) circle (3pt);
\end{tikzpicture}
\hspace*{5mm}
\begin{tikzpicture}[scale=.6]
\draw[->,line width=1pt] (0,0) -- (10,0);
\draw[->,line width=1.5pt] (0,0) -- (0,6);
\draw[line width=1pt] (0,2) .. controls (2,5) and (3,7) .. (5,3)
            .. controls (5,3) and (8,-3) .. (10,3);
\foreach \x in {.5,1,1.5,2,2.5,3,3.5,4,4.5,5,5.5,6,6.5,7,7.5,8,8.5,9,9.5}
  \draw[dashed] (\x,0) -- +(0,6);
\draw[line width=1pt] (0,0) -- (10,0);
\draw[line width=1pt] (0,1) -- (10,1);
\draw[line width=1pt] (0,2) -- (10,2);
\draw[line width=1pt] (0,3) -- (10,3);
\draw[line width=1pt] (0,4) -- (10,4);
\draw[line width=1pt] (0,5) -- (10,5);
\draw(.5,2.75) circle (3pt);
\draw(1,3.5) circle (3pt);   \filldraw(1,3) circle (3pt);
\draw(1.5,4.2) circle (3pt); \filldraw(1.5,4) circle (3pt);
\draw(2,4.75) circle (3pt);
\draw(2.5,5.1) circle (3pt);
\draw(3,5.2) circle (3pt);
\draw(3.5,5.05) circle (3pt);
\draw(4,4.6) circle (3pt);   \filldraw(4,5) circle (3pt);
\draw(4.5,3.9) circle (3pt); \filldraw(4.5,4) circle (3pt);
\draw(5,3) circle (3pt);     \filldraw(5,3) circle (3pt);
\draw(5.5,2.1) circle (3pt);
\draw(6,1.45) circle (3pt);  \filldraw(6,2) circle (3pt);
\draw(6.5,.95) circle (3pt);
\draw(7,.6) circle (3pt);
\draw(7.5,.35) circle (3pt);
\draw(8,.35) circle (3pt);
\draw(8.5,.5) circle (3pt);
\draw(9,.95) circle (3pt);
\draw(9.5,1.8) circle (3pt); \filldraw(9.5,1) circle (3pt);
\end{tikzpicture}
\end{tabular}
\end{center} \vspace{-1.5cm}
\caption{\label{fig:AADC}Non uniform sampling with an AADC \cite{Allier-Sicard-Fesquet-Renaudin05}.
Left: the disks correspond to the samples and the circle to the values at the clock ticks. Right: the disks are the only samples that are kept after decimation.}
\end{figure}
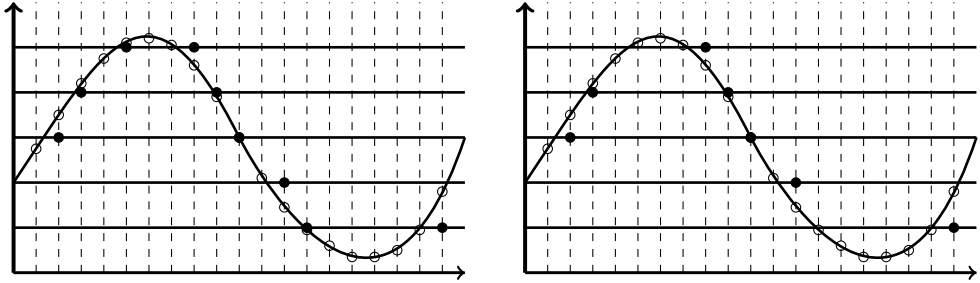

We introduce here an algorithm which is slightly different from the usual level crossing technique.
In \cite{Allier-Sicard-Fesquet-Renaudin05} the amplitudes are selected thanks a $M$-bit asynchronous analog-to-digital converter (AADC) that corresponds to $2^M$ predefined levels in the voltage range.
Level crossing sampling consists in taking a sample each time the predefined levels are crossed. Each amplitude has to be associated to a time. More precisely we store a delay elapsed since the last sample was taken, the local clock that enables this is then reset to zero and ready to measure the next delay.
In Figure \ref{fig:AADC} we display the case when the captured time is that of the next clock tick. The sample is displayed with disks, and the value of the signal at the capture times with circles (Figure \ref{fig:AADC}, left). This leads to a few number of samples, especially for sporadic signals. This procedure is refined decimating the samples by keeping only the last one when a level has been crossed many times successively (Figure \ref{fig:AADC}, right).

Our goal here is not to study the approximation of the signal but the number of non uniform samples, given the regularity of the signal, the clock precision and the level quantum. We introduce another sampling algorithm which is slightly different from the AADC, but easier to analyze mathematically, and which yields essentially the same number of samples. In Section \ref{sec:algo}, we describe this sampling algorithm and rephrase it in mathematical terms. In Section \ref{sec:Holder} we define the functions that we will use in the numerical experiments of Section 4. These functions are chosen because we are able to control exactly their H\"older regularity.

\section{Algorithm and mathematical interpretation}
\label{sec:algo}

\subsection{Step 1: Generation of an oversampled signal}
\label{sec:step1}

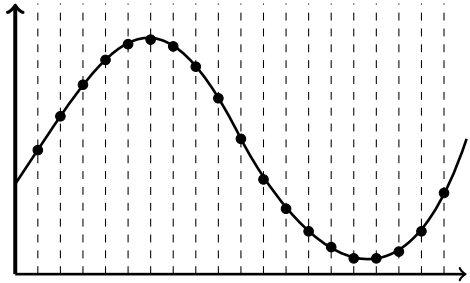
\begin{figure}[h]
\begin{center}
\begin{tikzpicture}[scale=.6]
\draw[->,line width=1pt] (0,0) -- (10,0);
\draw[->,line width=1.5pt] (0,0) -- (0,6);
\draw[line width=1pt] (0,2) .. controls (2,5) and (3,7) .. (5,3)
            .. controls (5,3) and (8,-3) .. (10,3);
\foreach \x in {.5,1,1.5,2,2.5,3,3.5,4,4.5,5,5.5,6,6.5,7,7.5,8,8.5,9,9.5}
  \draw[dashed] (\x,0) -- +(0,6);
\filldraw(.5,2.75) circle (3pt);
\filldraw(1,3.5) circle (3pt);
\filldraw(1.5,4.2) circle (3pt);
\filldraw(2,4.75) circle (3pt);
\filldraw(2.5,5.1) circle (3pt);
\filldraw(3,5.2) circle (3pt);
\filldraw(3.5,5.05) circle (3pt);
\filldraw(4,4.6) circle (3pt);
\filldraw(4.5,3.9) circle (3pt);
\filldraw(5,3) circle (3pt);
\filldraw(5.5,2.1) circle (3pt);
\filldraw(6,1.45) circle (3pt);
\filldraw(6.5,.95) circle (3pt);
\filldraw(7,.6) circle (3pt);
\filldraw(7.5,.35) circle (3pt);
\filldraw(8,.35) circle (3pt);
\filldraw(8.5,.5) circle (3pt);
\filldraw(9,.95) circle (3pt);
\filldraw(9.5,1.8) circle (3pt);
\end{tikzpicture}
\end{center} \vspace{-1.5cm}
\caption{\label{fig:step1}Regular sampling of the input signal (Step 1).}
\end{figure}

Even in event-driven systems, where the signal is not sampled at each clock tick, there are clocks that measure time, and specifically the time elapsed since the last event. These clocks have a certain precision, and all measured times are multiples of some basis time $t_b$.
Up to some re-scaling of time we suppose that $t_b=2^{-j}$, for some $j\in\mathbb{N}$.
We never have a complete knowledge of the original signal $f(t)$, but only its samples $f_{j,k}=f(k2^{-j})$, for all $k\in\bbZ$ (Figure \ref{fig:step1}).

Let $V_j$ be the space of continuous functions, which are linear on intervals
\begin{equation*}
I_{j,k} = [k2^{-j},(k+1)2^{-j}[ \text{ for all } k\in\bbZ.
\end{equation*}
The Faber--Schauder hierarchical basis, defined in \cite{Coh00}, yields a natural basis of $V_j$.
Let $\ph(x)=\max\{0,1-|x|\}$.
The functions $\ph_{j,k}=\ph(2^j\cdot-k)$, for all $k\in\bbZ$, form the Faber--Schauder basis.
We can uniquely define the linear interpolation $f_j\in V_j$ of $f$ at scale $2^{-j}$ by imposing $f_j(k2^{-j})=f_{j,k}$, for all $k\in\bbZ$, and
\begin{equation}
\label{eq:FS}
f_j =\sum_{k\in\bbZ} f_{j,k} \ph_{j,k}.
\end{equation}
In the sequel we suppose that $f$ is compactly supported in $[0,1]$ and therefore $k=0,\dots,2^j-1$ in Equation \eqref{eq:FS}.

\subsection{Step 2: Level crossing}
\label{sec:step2}

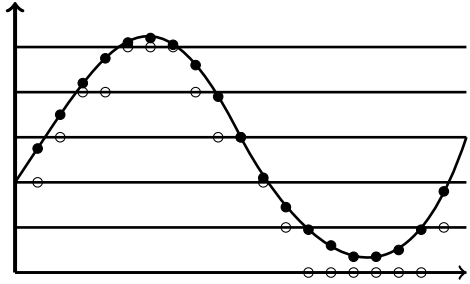
\begin{figure}[h]
\begin{center}
\begin{tikzpicture}[scale=.6]
\draw[->,line width=1pt] (0,0) -- (10,0);
\draw[->,line width=1.5pt] (0,0) -- (0,6);
\draw[line width=1pt] (0,2) .. controls (2,5) and (3,7) .. (5,3)
            .. controls (5,3) and (8,-3) .. (10,3);
\draw[line width=1pt] (0,0) -- (10,0);
\draw[line width=1pt] (0,1) -- (10,1);
\draw[line width=1pt] (0,2) -- (10,2);
\draw[line width=1pt] (0,3) -- (10,3);
\draw[line width=1pt] (0,4) -- (10,4);
\draw[line width=1pt] (0,5) -- (10,5);
\filldraw(.5,2.75) circle (3pt);  \draw(.5,2) circle (3pt);
\filldraw(1,3.5) circle (3pt);    \draw(1,3) circle (3pt);
\filldraw(1.5,4.2) circle (3pt);  \draw(1.5,4) circle (3pt);
\filldraw(2,4.75) circle (3pt);   \draw(2,4) circle (3pt);
\filldraw(2.5,5.1) circle (3pt); \draw(2.5,5) circle (3pt);
\filldraw(3,5.2) circle (3pt);    \draw(3,5) circle (3pt);
\filldraw(3.5,5.05) circle (3pt); \draw(3.5,5) circle (3pt);
\filldraw(4,4.6) circle (3pt);    \draw(4,4) circle (3pt);
\filldraw(4.5,3.9) circle (3pt);  \draw(4.5,3) circle (3pt);
\filldraw(5,3) circle (3pt);      \draw(5,3) circle (3pt);
\filldraw(5.5,2.1) circle (3pt);  \draw(5.5,2) circle (3pt);
\filldraw(6,1.45) circle (3pt);   \draw(6,1) circle (3pt);
\filldraw(6.5,.95) circle (3pt);  \draw(6.5,0) circle (3pt);
\filldraw(7,.6) circle (3pt);     \draw(7,0) circle (3pt);
\filldraw(7.5,.35) circle (3pt);  \draw(7.5,0) circle (3pt);
\filldraw(8,.35) circle (3pt);    \draw(8,0) circle (3pt);
\filldraw(8.5,.5) circle (3pt);   \draw(8.5,0) circle (3pt);
\filldraw(9,.95) circle (3pt);    \draw(9,0) circle (3pt);
\filldraw(9.5,1.8) circle (3pt);  \draw(9.5,1) circle (3pt);
\end{tikzpicture}
\end{center} \vspace{-1.5cm}
\caption{\label{fig:step2}Reduction to predefined levels (Step 2). Samples from Step 1 are disks and the new samples are the circles.}
\end{figure}

We consider that levels are uniformly spaced by some quantum $2^{-M}$. In applications where the range of the signal is $[0,1]$, the sample can then be stored with a $M$-bit register.
The second step consists in approximating the samples $f_{j,k}$ by the nearest level below (Figure \ref{fig:step2}).

We denote $\floor{x}$ the integer part of $x$, namely $\floor{x}=\sup\{n\in\bbN,\,n\leq x\}$.
The function $\tilde f_j\in V_j$ which coincides with the new samples is uniquely defined by
\begin{equation*}
\tilde f_j = \sum_{k=0}^{2^j-1} 2^{-M} \floor{2^M f(k2^{-j})} \ph_{j,k}.
\end{equation*}

\subsection{Step 3: Decimation}

\begin{figure}[h]
\begin{center}
\begin{tikzpicture}[scale=.6]
\draw[->,line width=1pt] (0,0) -- (10,0);
\draw[->,line width=1.5pt] (0,0) -- (0,6);
\draw[line width=1pt] (0,2) .. controls (2,5) and (3,7) .. (5,3)
            .. controls (5,3) and (8,-3) .. (10,3);
\draw[line width=1pt] (0,0) -- (10,0);
\draw[line width=1pt] (0,1) -- (10,1);
\draw[line width=1pt] (0,2) -- (10,2);
\draw[line width=1pt] (0,3) -- (10,3);
\draw[line width=1pt] (0,4) -- (10,4);
\draw[line width=1pt] (0,5) -- (10,5);
\filldraw(.5,2) circle (3pt);
\filldraw(1,3) circle (3pt);
\draw(1.5,4) circle (3pt);
\filldraw(2,4) circle (3pt);
\draw(2.5,5) circle (3pt);
\draw(3,5) circle (3pt);
\filldraw(3.5,5) circle (3pt);
\filldraw(4,4) circle (3pt);
\draw(4.5,3) circle (3pt);
\filldraw(5,3) circle (3pt);
\filldraw(5.5,2) circle (3pt);
\filldraw(6,1) circle (3pt);
\draw(6.5,0) circle (3pt);
\draw(7,0) circle (3pt);
\draw(7.5,0) circle (3pt);
\draw(8,0) circle (3pt);
\draw(8.5,0) circle (3pt);
\filldraw(9,0) circle (3pt);
\filldraw(9.5,1) circle (3pt);
\end{tikzpicture}
\end{center} \vspace{-1.5cm}
\caption{\label{fig:step3}Non uniform samples after decimation (Step 3).}
\end{figure}
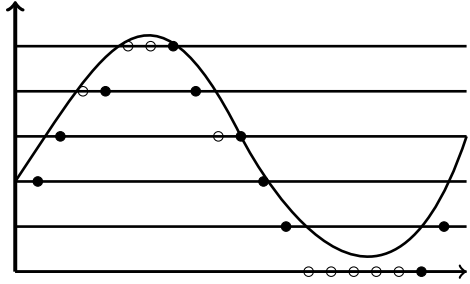

Next, we decimate the samples so as to keep only one sample when consecutive samples have the same amplitude. We choose to keep the last sample to be compatible with the causality principle (Figure \ref{fig:step3}).

We are interested in the number of samples after the three steps. Comparing Figures \ref{fig:AADC} (right) and \ref{fig:step3}, we notice that the number of final samples (the disks in both figures) are comparable. This is a generic situation. In fact the differences are mainly due to the extreme upper and lower levels.

From the mathematical point of view, decimation consists in keeping a subsequence of $k=0,\dots,2^j-1$, defined by induction: $k_0=0$ and
\begin{equation*}
k_{i+1}=\min \{k\geq 1+k_i\slash\floor{2^M f(k2^{-j})}\neq \floor{2^M f(k_i2^{-j})}\}.
\end{equation*}
We only store the couples $(\delta t_i,a_i)$ where $\delta t_i=(k_i-k_{i-1})2^{-j}$, $i\geq 1$ is the delay since the last sample, and $a_i=2^{-M} \floor{2^M f(k_i2^{-j})})$ is the amplitude of the sample.
Step 3 leads to a reduction of the number of samples, but does not introduce any approximation.

\section{Application to monoH\"olderian functions}
\label{sec:Holder}

Our goal is to relate the number of non uniform samples to the regularity of the signal. In particular, we address the H\"olderian regularity.

\subsection{MonoH\"olderian functions}

Before introducing the H\"olderian regularity, we first recall a few definitions. They allow a weaker definition of pointwise smoothness. The final goal is to define strongly monoH\"olderian functions, a notion that formalizes the idea of a function which has the as uniformly as possible regularity.

Let $x_0\in\bbR$ and $r>0$, for all $0\leq h\leq r$ we define
\begin{equation*}
B_h(x_0,r)=\{x:[x,x+h]\subset B(x_0,r)\},
\end{equation*}
where $B(x_0,r)$ is the ball of center $x_0$ and radius $r$, and denote, as usual, the oscillation of a function $f:\bbR\to\bbR$ at $x_0$ on the ball $B(x_0,r)$ as
\begin{equation*}
\Osc(f)(x_0,r)=\sup_{|h|\leq r} \| f(x+h)-f(x) \|_{L^\infty(B_h(x_0,r))}.
\end{equation*}

\begin{definition}
Let $f:\bbR\to\bbR$ be a locally bounded function, let $x_0\in\bbR$ and $\alpha\in(0,1)$. The function $f$ is H\"olderian of exponent $\alpha$ at $x_0$ ($f\in C^\alpha(x_0)$) if there exist $C$ and $R>0$ such that
\begin{equation}
\label{eq:strong holder-fd}
\Osc(f)(x_0,r)\leq C r^\alpha, \quad \forall r\leq R.
\end{equation}
A function $f$ is uniformly H\"olderian of exponent $\alpha$ ($f\in C^\alpha(\bbR)$) if $C$ and $R$ in Equation (\ref{eq:strong holder-fd}) are uniform in $x_0\in\bbR$.
\end{definition}

The irregularity of a function can be studied through the notion of anti-H\"olderianity.

\begin{definition}
Let $f:\bbR\to\bbR$ be a locally bounded function, let $x_0\in\bbR$ and $\alpha\in(0,1)$. The function $f$ is anti-H\"olderian of exponent $\alpha$ at $x_0$ ($f\in I^\alpha(x_0)$) if there exist $C$ and $R>0$ such that
\begin{equation}
\label{eq:anti-holder}
\Osc(f)(x_0,r)\geq C r^\alpha, \quad \forall r\leq R.
\end{equation}
\end{definition}

Let us notice that the statement (\ref{eq:anti-holder}) is stronger than just negating the H\"olderian regularity. Indeed such a negation only yields the existence, for any $C>0$, of a subsequence $(r_n)_n$ (depending on $C$) for which
\begin{equation*}
\Osc(f)(x_0,r_n)\geq C r_n^\alpha.
\end{equation*}
Strongly monoH\"olderian functions naturally arise in the study of the regularity of mappings such as Weierstrass-type or random processes (see e.g.\ \cite{gem-hor:80,heurt:03}). Indeed, many results only hold for such mappings.

\begin{definition}
Let $\alpha\in(0,1)$. A function $f:\bbR\to\bbR$ is strongly monoH\"olderian of exponent $\alpha$ ($f\in SM^\alpha(\bbR)$) if $f\in C^\alpha(\bbR)\cap I^\alpha(\bbR)$, i.e.\ if there exists $C$ and $R>0$ such that, for any $x_0\in\bbR$,
\begin{equation}\label{e:monoholder}
 r^\alpha/C \le \Osc(f)(x_0,r)\le Cr^\alpha \quad \forall r\le R.
\end{equation}
\end{definition}

\subsection{Approximation properties}

As already mentioned, only Steps 1 and 2 lead to approximations. To state our approximation results, we need some preliminary definitions. According to our application, we now restrict to functions defined on $[0,1]$. For any continuous function $f$ on $[0,1]$, we define its uniform regularity modulus by
\begin{equation*}
\omega_f(r)=\sup_{|h|\leq r,\,x\slash[x,x+h]\subset (0,1)}|f(x+h)-f(x)|.
\end{equation*}
The function $\omega_f$ is {\it a modulus of continuity} in the sense that
$\omega_f(0)=0$ and that there exists some $C>0$ such that $\omega_f(2r)\leq C\omega_f(r)$ (see~\cite{JafMey96}).

In what follows we need the notion of {\it strong modulus of continuity} introduced in~\cite{Clau07,JafMey96}. The modulus of continuity $\theta$ is said to be strong if there exists $C>0$ such that for any positive integer $J$ one has
\begin{equation*}
\sum_{j=0}^J 2^j\theta(2^{-j})\leq C2^J\theta(2^{-J})
\text{ and }
\sum_{j=J}^\infty \theta(2^{-j})\leq C\theta(2^{-J}).
\end{equation*}
It is well-known \cite{JafMey96} that if there exists some strong modulus of continuity $\theta$ such that
\begin{equation*}
\omega_f(2^{-j})\leq C\theta(2^{-j}),
\end{equation*}
then
\begin{equation*}
\|f-f_j\|_{L^\infty} \leq C \theta(2^{-j}),
\end{equation*}
where $f_j$ is defined by Equation \eqref{eq:FS}.
In particular if $f\in C^\alpha(0,1)$ there exists some $C>0$ such that for any $j\geq 0$
\begin{equation*}
\omega_f(2^{-j})\leq C\theta(2^{-j}),
\end{equation*}
with $\theta(2^{-j})=2^{-j\alpha}$, and then there exists a constant $C$ (which depends on $f$ but not on the scale $j$) such that
\begin{equation*}
\|f-f_j\|_{L^\infty} \leq C 2^{-j\alpha}.
\end{equation*}
Assume now that in addition
\begin{equation*}
\omega_f(2^{-j})\geq \theta(2^{-j})/C,
\end{equation*}
then, following \cite{Clau07,CN12}, there exists a constant $C$ and $\beta>1$ (which depend on $f$ but not on the scale $j$) such that\begin{equation*}
\|f-f_j\|_{L^\infty} \geq j^{-\beta}\theta(2^{-j})/C.
\end{equation*}
In particular, applying these results with the strong modulus of continuity $\theta(r)=r^\alpha$, we deduce that if the function $f$ is assumed to be uniformly monoH\"olderian with exponent $\alpha$ there exist some $C>0$ and $\beta>1$ such that
\begin{equation*}
2^{-j\alpha}j^{-\beta}/C\leq \|f-f_j\|_{L^\infty} \leq C2^{-j\alpha}.
\end{equation*}
This yields estimates on the error due to Step 1. The approximation made at Step 2 clearly does not depend on the regularity of function $f$, and we have
\begin{equation*}
\|f_j-\tilde f_j\|_{L^\infty} \leq 2^{-M}.
\end{equation*}

\subsection{Theoretical number of samples in the case of a monotonous function}

If $f$ is a monoH\"olderian function with exponent $\alpha$, by definition there exists $C_1,C_2>0$ and for any scale $j\geq 0$ and $0\leq k\leq 2^j-1$
\begin{eqnarray*}
C_1 2^{-j\alpha}\leq \sup_{(x,h)\slash[x,x+h]\subset\left[\frac k{2^j},\frac{k+1}{2^j}\right]}|f(x+h)-f(x)|\leq C_2 2^{-j\alpha}.
\end{eqnarray*}
If the function is additionally supposed to be monotonous, we further have exactly
\begin{equation*}
\sup_{(x,h)\slash[x,x+h]\subset\left[\frac k{2^j},\frac{k+1}{2^j}\right]}|f(x+h)-f(x)|
=\left|f\left(\frac{k+1}{2^j}\right)-f\left(\frac{k}{2^j}\right)\right|.
\end{equation*}
Hence
\begin{equation*}
C_1 2^{j(1-\alpha)}
\leq |f(1)-f(0)|
= \sum_{k=0}^{2^{j}-1} \left|f\left(\frac{k+1}{2^j}\right)-f\left(\frac{k}{2^j}\right)\right|
\leq C_22^{j(1-\alpha)}.
\end{equation*}
A monoH\"olderian signal crosses equi-spaced levels with quantum $2^{-M}$ at most $C_2 2^{M+(1-\alpha)j}$ times. The worst case is that of monotonous signals.

Besides, initial sampling (Step 1) takes exactly $2^j$ samples. This is hence the first natural upper bound for the number of samples. Together with monoH\"older\-ianity we know that the number of samples is less than the minimum of these two bounds.
For large values of $M$ (or small values of $\alpha$), we indeed keep almost all of the $2^j$ samples. Otherwise we can expect some reduction of the number of samples.
For $C=1$, the threshold is $M\simeq\alpha j$. Observe that the proof is based on the fact, that in the monotonous case, we can estimate in a very simple way the oscillations
\begin{equation*}
\sup_{(x,h)\slash[x,x+h]\subset\left[\frac k{2^j},\frac{k+1}{2^j}\right]}|f(x+h)-f(x)|
\end{equation*}
of the function. Of course in the general case, the situation can be much more complicated. Nevertheless, generic results in the sense of prevalence as stated in \cite{CN10} are expected to hold. In what follows, we illustrate through numerical simulations what happens in two cases.

\section{Numerical simulations}

\subsection{Fractional Brownian motion and the Weierstrass function}

We test level crossing on two toy models: sample paths of fractional Brownian motion $B_H$ and the Weierstrass function $W_H$, which are indexed by the Hurst index $H\in(0,1)$. The choice of these two cases is guided by the fact that their smoothness properties are related to the Hurst index.

\begin{figure}[h]
\includegraphics[width=\textwidth,height=4cm]{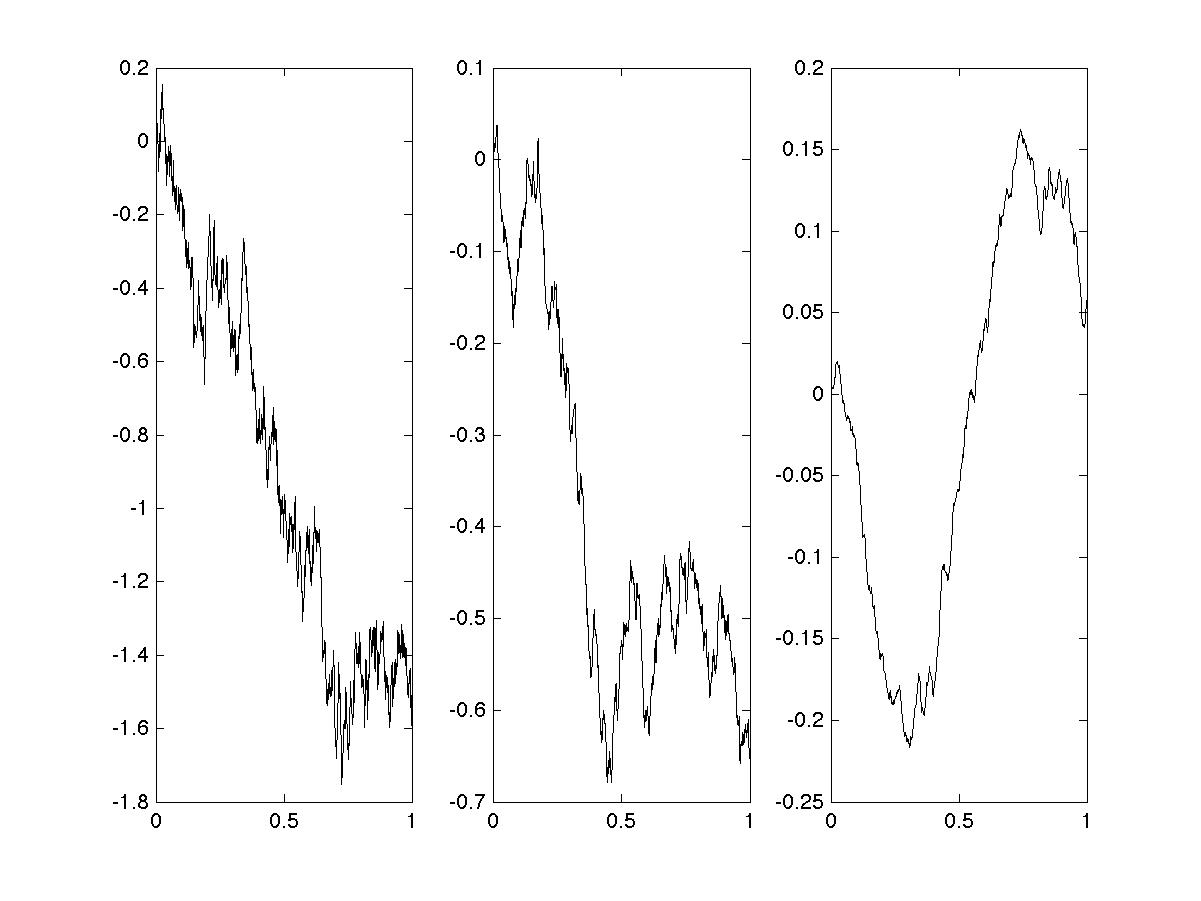}
\caption{\label{fig:fBm}Three realizations of fractional Brownian motions for $H=0.5$, $H=0.7$, and $H=0.9$ (from left to right).}
\end{figure}

The \textbf{fractional Brownian motion} (fBm) $B_H$ is the unique Gaussian $H$-self-similar process with stationary increments. It is defined from its covariance function
\begin{equation*}
\bbE\big[B_H(x)B_H(y)\big]=\tfrac{1}{2} \left(|x|^{2H}+|y|^{2H}-|x-y|^{2H}\right)
\end{equation*}
for all $(x,y)\in[0,1]^2$. The classical Brownian motion corresponds to $H=1/2$.
The sample paths of fBm are well-known to be almost surely continuous. Further, its Hurst index $H$ is directly related to the roughness of its sample paths. More precisely the classical law of the iterated logarithm ensures that
\begin{equation*}
B_H\in C^{H-\varepsilon}([0,1])\cap I^{H+\varepsilon}([0,1]) \text{ almost surely}.
\end{equation*}
Roughly speaking, a.s. for all $(x,y)\in[0,1]^2$,
\begin{equation*}
\sup_{(u,v)\in [x,y]^2}|B_H(u)-B_H(v)|\sim|x-y|^H.
\end{equation*}
Figure \ref{fig:fBm} presents three realizations of sample paths of fractional Brownian motions for $H=0.5$, $H=0.7$, $H=0.9$, and $j=10$ ($1024$ samples).

\begin{figure}[h]
\includegraphics[width=\textwidth,height=4cm]{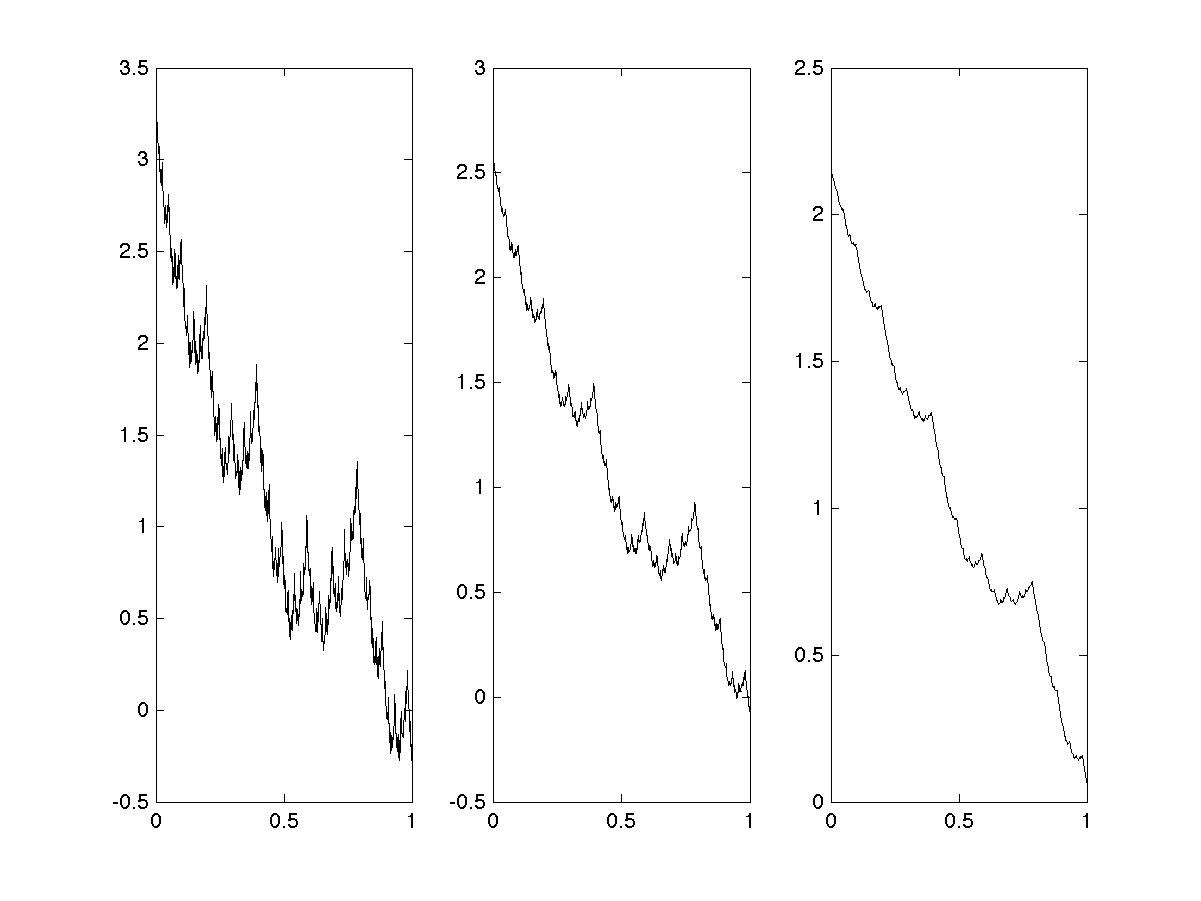}
\caption{\label{fig:Weierstrass}The Weierstrass function for $H=0.5$, $H=0.7$, and $H=0.9$ (from left to right).}
\end{figure}

The \textbf{Weierstrass function} $W_H$ is a classical example of monoH\"olderian function with exponent $H$ as proved in~\cite{Har16}. It is defined as
\begin{equation*}
W_H(x) = \sum_{j=0}^\infty 2^{-jH} \cos(2^jx),
\end{equation*}
and, for all $(x,y)\in[0,1]^2$,
\begin{equation*}
\sup_{(u,v)\in [x,y]^2}|W_H(u)-W_H(v)|\sim|x-y|^H.
\end{equation*}
Figure \ref{fig:Weierstrass} presents the graphs of the Weierstrass functions for $H=0.5$, $H=0.7$, $H=0.9$, and $j=10$ ($1024$ samples).

\subsection{Tests}

The tests are performed within the SPASS \textsc{Matlab} toolbox \cite{SPASS} (Signal Processing for ASynchronous Systems toolbox). It has been originally designed to treat non uniform signals produced by asynchronous systems, but can be used for a large variety of signals.
To generate fractional Brownian motions, we make use of the \texttt{genFBMJFC.m} function \cite{Jeff}.

We use two values of $j$ (10 and 13) and two values of $M$ (4 and 5).
These small values of $M$ are sufficient for most mobile applications.
Our output is the number of samples after decimation (Step 3).
For the fractional Brownian motion, we perform 1000 realizations and average the number of samples obtained for each realization to obtain an average number $n$. We perform the same tests on the Weierstrass function (deterministic function, only one realization).

\begin{figure}[h]
\includegraphics[width=\columnwidth,height=4cm]{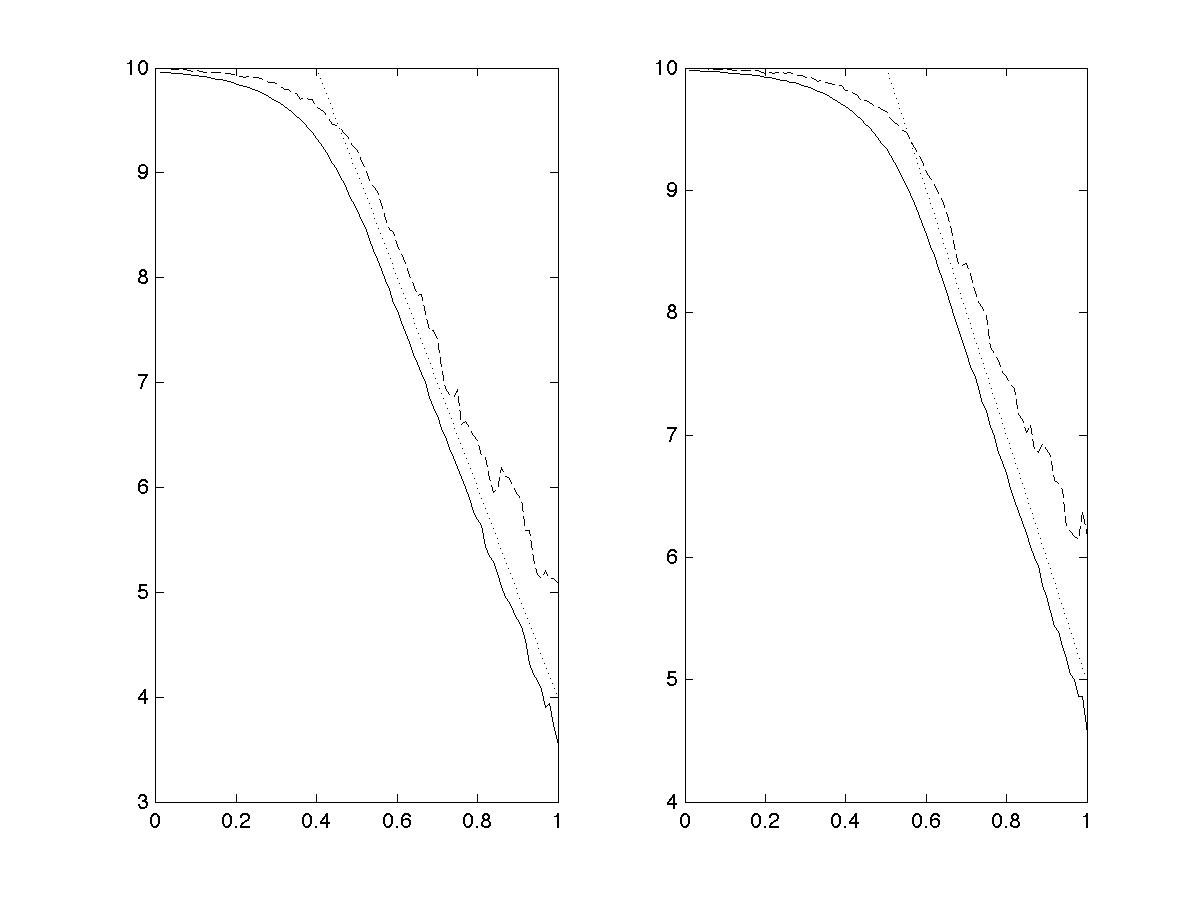}
\caption{\label{fig:complexity_10}Number $n$ of samples in terms of the Hurst number in the log scale for $j=10$, and $M=4$ (left) and $M=5$ (right). Solid lines correspond to the averaged number for the fractional Brownian motion, the dashed lines to the Weierstrass function, and the dotted lines to the worst case $j+\max(1,M-Hj)$.}
\end{figure}

\begin{figure}[h]
\includegraphics[width=\columnwidth,height=4cm]{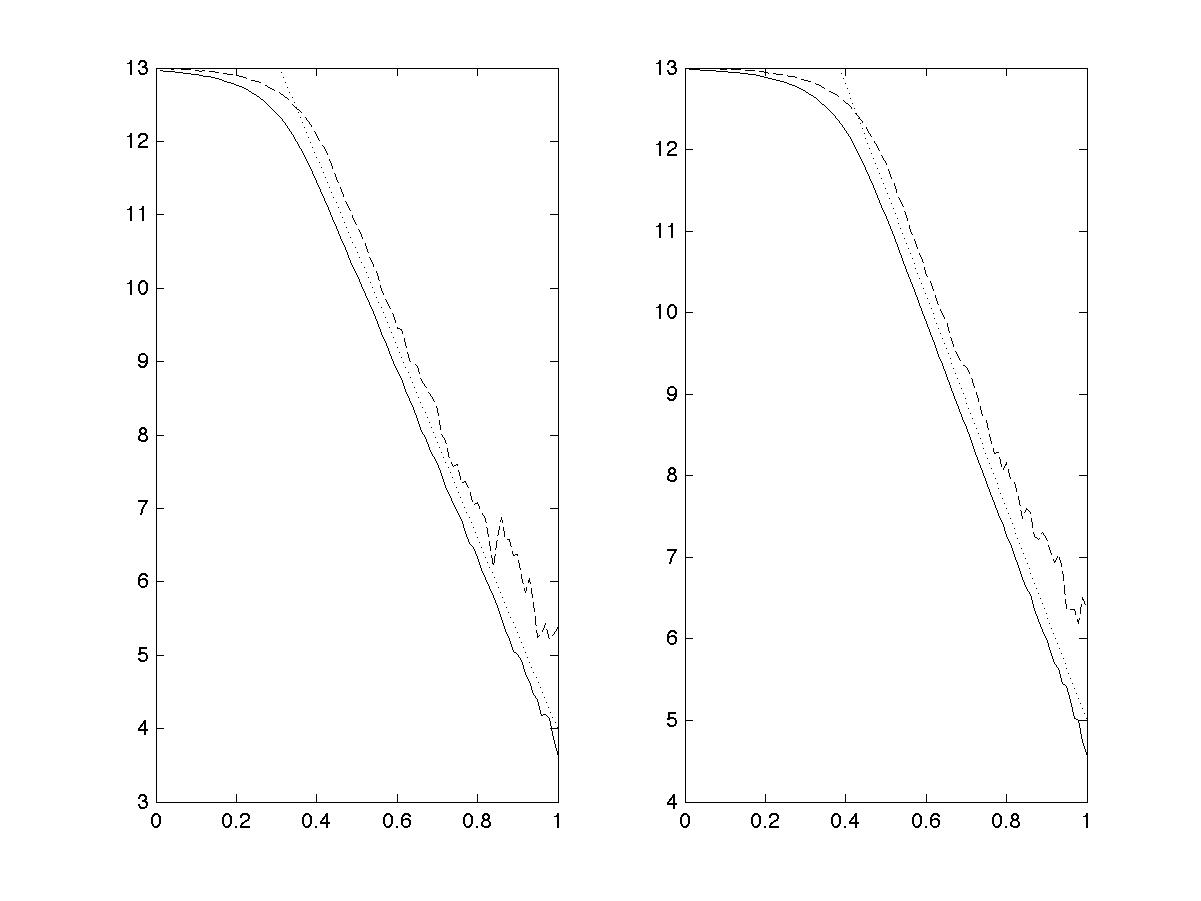}
\caption{\label{fig:complexity_13}Number $n$ of samples in terms of the Hurst number in the log scale for $j=13$, and $M=4$ (left) and $M=5$ (right). Same plotting conventions as Figure \ref{fig:complexity_10}.}
\end{figure}

We perform this for values of the Hurst number $H$ in the $(0,1)$ range and obtain the plots (in semi-log scale, with log-basis 2) in Figures \ref{fig:complexity_10} and \ref{fig:complexity_13} for $j=10$ and 13 respectively.
We also plot the number of samples computed in the worst case (monotonous function i.e. maximum total variation) for $C=1$: $N_{\rm worst}=\min(2^j,2^{M+(1-H)j})=2^j\min(1,2^{M-Hj})$.

We distinguish two regimes: below some value of the Hurst number $H\sim M/j$ the algorithm more or less keeps all the original samples, above this value the decimation is efficient and yields a significant reduction of the number of samples.
For the different curves these "critical" values of $H$ are given in Table \ref{tab:critical}.

\begin{table}[h]
\begin{center}
\begin{tabular}{c|c|c}
& $M=4$ & $M=5$ \\
\hline
$j=10$ & 0.4 & 0.5 \\
\hline
$j=13$ & $\sim 0.3$ & $\sim 0.4$
\end{tabular}
\end{center}
\caption{\label{tab:critical}"Critical" values of the Hurst number.}
\end{table}

The plots associated to fBm are much more regular than those associated to the Weierstrass function because they are obtained by an averaging procedure. Besides, for the Weierstrass function, the constant involved in Equation~\eqref{e:monoholder} is \textit{a priori} not equal to $1$ and indeed depends on $H$. It can be explained using the range of the fractional derivative of order $H$ of Weierstrass function $W_H$ which is not reduced to a constant and depends on $H$ (\cite{KWS02,SKM93,ZZ96} for more details). 
 
\section{Conclusion}

We have predicted for monoH\"olderian functions and shown numerically that there is a strong relationship between the smoothness properties of a signal and the number of samples that can be obtained by the crossing level algorithm presented in this paper.
This is rigorously proved in the case of monotonous monoH\"olderian functions. The next step, which will be the purpose of a forthcoming paper, will then be to consider signals whose regularity may change from point to point such as multifractional or multifractal signals.

\vspace{13pt}
\centerline{ACKNOWLEDGEMENT}
\vspace{13pt}

\noindent LJK is partner of the LabEx PERSYVAL-Lab (ANR--11-LABX-0025-01) funded by the French program Investissement d'avenir and this work has been partially supported by the MathSTIC Project OASIS of the Grenoble University. The authors wish to thank Jean-Fran\c{c}ois Coeurjolly for fruitful discussions.

\end{document}